\documentclass[11pt]{amsart}
\usepackage{amsmath,amscd,amssymb,amsfonts,amsthm,hyperref}

\usepackage[cmtip,matrix,arrow]{xy}
\usepackage{tikz-cd}
\usepackage{tikz}
\usetikzlibrary{decorations.markings}
\usetikzlibrary{arrows}

\newtheorem{theorem}{Theorem}[section]

\theoremstyle{definition}    
\newtheorem{definition}[theorem]{Definition}
\theoremstyle{remark}

\renewcommand{\L}{\mathcal{L}}
\renewcommand{\O}{\mathcal{O}}

\newcommand{\ca}{\mathcal}

\newcommand{\E}{\ca{E}}

\newcommand{\F}{\mathcal{F}}

\newcommand{\R}{\mathbb{R}}

\newcommand{\SU}{\on{SU}}

\newcommand{\Z}{\mathbb{Z}}

\newcommand\pt{\on{pt}}

\newcommand{\zz}{\phantom{\cdot}^0 }
\newcommand{\bb}{\phantom{\cdot}^b }
\newcommand\lie[1]{\mathfrak{#1}}
\renewcommand{\k}{\lie{k}}
\newcommand{\h}{\lie{h}}
\newcommand{\g}{\lie{g}}

\renewcommand{\a}{\mathsf{a}}

\newcommand{\on}{\operatorname}

\renewcommand{\ker}{ \on{ker}}

\newcommand{\SO}{ \on{SO}}

\newcommand{\sz}{\mathsf{s}}
\newcommand{\tz}{\mathsf{t}}

\newcommand{\da}{\dasharrow}

\newcommand\qu{/\kern-.7ex/} 

\renewcommand{\d}{{\mbox{d}}}

\newcommand{\f}{\frac}

\newcommand\hh{{\f{1}{2}}}

\newcommand{\eeq}{\end{eqnarray*}}
\newcommand{\beq}{\begin{eqnarray*}}

\newcommand{\pr}{\on{pr}}
\newcommand{\wh}{\widehat}

\newcommand{\mf}{\mathfrak}
\newcommand{\rra}{\rightrightarrows}
\renewcommand{\subset}{\subseteq}

\newcommand{\wt}{\widetilde}
\newcommand{\Ra}{\Rightarrow}
\newcommand{\VB}{{\mathcal{VB}}}

\newcommand{\LA}{\mathcal{LA}}

\begin{document}

\title{Lie algebroids}

\author{Eckhard Meinrenken}
\address{Mathematics Department\\ University of Toronto\\40 St George Street\\Toronto, ON M5S2E4}
\email{mein@math.toronto.edu}
\begin{abstract} This is an overview article on Lie algebroids, and their role as the infinitesimal counterparts of Lie groupoids. 
\end{abstract}
\maketitle

\section{Introduction}
Lie algebroids are differential-geometric structures unifying tangent bundles and Lie algebras.  They appear naturally in foliation theory, Poisson geometry, geometric mechanics, generalized complex geometry, and many other fields of mathematics and physics. 
Lie algebroids give rise to infinitesimal symmetries of manifolds, generalizing actions of Lie algebras.
Their global counterparts are Lie groupoids, and one of the major themes in the theory is the interplay between the infinitesimal and global aspects. Throughout the following overview, we will work with (real) $C^\infty$-manifolds; similar notions exist in the complex 
holomorphic setting.   

\section{Definitions and basic properties} 
\subsection{Definition and examples}\label{subsec:defex}
Let $M$ be a manifold, and denote by $\L_X$ the Lie derivative with respect to a vector field $X$. 
 \begin{definition}[Pradines \cite{pra:th}]
 	A \emph{Lie algebroid} is a vector bundle $A\to M$, together with a 
 	Lie bracket $[\cdot,\cdot]$ on its space of sections and a 
 	vector bundle map 
 	\[ \a\colon A\to TM \]
 	called the \emph{anchor}, such that the  \emph{Leibniz rule}
 	\[ [\sigma,f\tau]=f[\sigma,\tau]+(\L_{\a(\sigma)}f)\ \tau\] holds 
 	for all $\sigma,\tau\in \Gamma(A)$ and all $f\in C^\infty(M)$. 
 \end{definition}
 We will use the notation $A\Rightarrow M$ for Lie algebroids; the double arrow is meant to suggest the 
 infinitesimal version of Lie groupoids $G\rra M$. Some standard examples:
 
 \begin{enumerate}
 	\item A Lie algebra is the same as a Lie 
 	algebroid $\g\Ra \pt$
 	over a point. More generally, Lie algebroids with zero anchor are \emph{families of Lie algebras}. (These need not be locally trivial, in 
 	general.) 
 	\item The tangent bundle is a Lie algebroid $TM\Ra M$ with anchor $\a=\on{id}_{TM}$. 
 	\item 
 	The subbundle of vectors tangent  to a foliation $\ca{F}$  is a Lie algebroid $T_{\ca{F}}M\Ra  M$, with anchor the inclusion. 
 	\item A principal $G$-bundle $P\to M$ defines an \emph{Atiyah algebroid} 
 	$$\on{At}(P)=(TP)/G\Ra M.$$
 	Its sections are the infinitesimal automorphisms of $P$, that is, 
  	$G$-invariant vector fields on $P$. 
 	For a vector bundle $V\to M$, one may define $\on{At}(V)\Ra M$ as the Atiyah algebroid of the frame 
 	bundle. Its sections are the infinitesimal automorphisms of $V$, that is, 
 	linear vector fields on the total space. Equivalently, they are the first order differential operators acting on sections of $V$, 
 	with scalar principal symbol. 
\item 
An action of a Lie algebra $\g$ on $M$ is a Lie algebra morphism $\rho\colon \g\to \mf{X}(M)$  such that 
$$\a\colon \g\times  M\to TM,\ (\xi,x)\mapsto \rho(\xi)_x$$ 
is smooth. It defines an \emph{action Lie algebroid} $$\g\times M\Ra M,$$ with the  $\a$ as the anchor, and with bracket on sections extending the given one on constant sections. 
 	\item\label{it:f} Given a hypersurface $Z\subset M$, there are Lie algebroids \cite{mel:ati} 
 	$$ \bb TM\Ra M,\ \ \zz TM\Ra M$$  whose sections are the vector fields  tangent to $Z$ and vanishing along $Z$, respectively. 
 	\item\label{it:g} A closed 2-form $\omega\in \Omega^2(M)$ defines an \emph{Almeida-Molino Lie algebroid} \cite{alm:sui1}
 	 $$A_\omega=TM\oplus \R,$$ with bracket on sections 
 	 $$ [X+f,Y+g]=[X,Y]+\L_X g-\L_Y f+\omega(X,Y).$$
 	\item A Poisson bracket $\{\cdot,\cdot\}$ on $M$ defines a \emph{cotangent Lie algebroid}
 	$$ T^*M\Ra M.$$
 	 The bracket on 1-forms is uniquely determined by its values on exact 1-forms 
 	$$ [\d f,\d g]=\d \{f,g\};$$
 	 the anchor takes $\d f$ to the Hamiltonian vector field $X_f=\{f,\cdot\}$. 
 \end{enumerate}	
%
The \emph{product}  of two Lie algebroids $A_1\Ra M,\ A_2\Ra M$ 
is uniquely a Lie algebroid $$A_1\times A_2\Ra  M_1\times M_2,$$ 
in such a way that 
the map $\Gamma(A_1)\oplus \Gamma(A_2)\to \Gamma(A_1\times A_2)$, $$
(\sigma_1,\sigma_2)\mapsto \pr_1^*\sigma_1+\pr_2^*\sigma_2$$
is a Lie algebra morphism. 

\subsection{Subalgebroids}

 A \emph{Lie subalgebroid} of  a Lie algebroid $A\Ra M$ is a vector subbundle $B\subseteq A$ along $N\subseteq M$, 
 such that:  
 \begin{itemize}
 	\item $\Gamma(A,B)=\{\sigma\in \Gamma(A)|\ \sigma|_N\in \Gamma(B)\}$ is a Lie subalgebra. 
 	\item The anchor of $A$ takes $B$ to $TN$.
 \end{itemize}
A Lie subalgebroid acquires a Lie algebroid structure $B\Rightarrow N$, with
the unique bracket such that restriction map $\Gamma(A,B)\to \Gamma(B),\ \sigma\mapsto \sigma|_N$ 
preserves brackets. 
 Examples:
 
 \begin{enumerate}
 	\item  If $A=\g$ is a Lie algebra, one recovers the usual notion of a Lie subalgebra. 
 	\item If $A\Ra M$ is any Lie algebroid, and $N\subset M$ is such that 
 	$B=\a^{-1}(TN)$  is a sub-bundle, then $B$ is a sub-Lie algebroid. In particular, 
  if  $\a\colon A\to TM$ has constant rank, then $\ker(\a)\subset A$ is a Lie subalgebroid.  
 	\item  Let $M$ be 	Poisson manifold. The conormal bundle of a submanifold $N\subset M$ 
 	is a Lie subalgebroid of $T^*M$ 
 	if and only if $N$ is a coisotropic submanifold.  
 \end{enumerate}

 \subsection{Morphisms and isomorphisms}

 Let $A\Rightarrow M,\ B\Rightarrow  N$ be Lie algebroids.  
 A \emph{Lie algebroid morphism} (respectively, \emph{comorphism})  from $B$ to $A$ along $f\colon N\to M$ is a 
Lie subalgebroid $$\Gamma\subset A\times B$$ along the graph $\on{Gr}(f)=\{(f(n),n)\}\subset M\times N$,
such that projection to $B$ (respectively, to $A$) is a fiberwise vector space isomorphism. 

Lie algebroid morphisms are, in particular, vector bundle morphisms $B\to A$; the description in terms of the graph is used to formulate compatibility with the bracket. A comorphism defines a bundle map  $f^*A\to B$; compatibility with the bracket means that the 
pullback map on sections $\Gamma(A)\to \Gamma(B)$ is a Lie algebra morphism. We use the notation 
$B\da A$ to indicate comorphisms; the direction of the arrow is given by the base map. 

	\begin{enumerate}
		\item 	Examples of Lie algebroid morphisms include: morphisms of Lie algebras, the tangent map $Tf\colon TN\to TM$ of $f\colon N\to M$, the inclusion of a Lie subalgebroid, the anchor map $A\to TM$ of any Lie algebroid. 
		\item A bundle map $\alpha\colon TM\to \g$ to a Lie algebra 
		is a Lie algebroid morphism if and only if the associated 1-form $\alpha\in \Omega^1(M,\g)$ is a Maurer-Cartan element, $\d \alpha+\hh [\alpha,\alpha]=0$. 	
	   \item If $P\to M$ is a principal bundle, a Lie algebroid morphism $TM\to \on{At}(P)$ (with base map the identity) is the same as a flat principal connection on $P$. 
		\item A Lie algebroid morphism $T[0,1]\to A$ from the tangent bundle of an interval to a Lie algebroid
		is called an \emph{$A$-path}. It is equivalent to a path $\gamma\colon [0,1]\to A$ with the property 
		\[ \a\circ \gamma=\f{d}{d t}(\pi\circ \gamma) \]
		where $\pi\colon A\to M$ is the bundle projection. 
		\item 	A  Lie algebroid comorphism $TM\da \g$ is the same as a $\g$-action  on $M$. Here $\Gamma\subset \g\times TM$ is identified with the action Lie algebroid. Generalizing, one defines an action of a Lie algebroid $A\Ra M$ on a manifold $N$ to be a 
		Lie algebroid comorphism $TN\da A$ along some map $f\colon N\to M$. 	
		\item Let $M_1,M_2$ be Poisson manifolds. 
				A Poisson map $M_1\to M_2$ gives a Lie algebroid comorphism $T^*M_1\da T^*M_2.$  
	\end{enumerate}
Lie algebroid morphisms can be composed, defining a category $\LA$. Similarly, comorphisms can be composed, defining a category $\LA^\vee$.

\subsection{Pullbacks}
Let $A\Ra M$ be a Lie algebroid, and $f\colon N\to M$ a smooth map transverse to the anchor of $A$. That 
is, for $x=f(y)$, 
$$T_xM=\on{ran}(T_yf)+\on{ran}(\a_x).$$  
Then the fiber product 
$ f^!A=A\times_{TM} TN$
is a vector bundle over $N$, and acquires 
a Lie algebroid structure $$ f^!A\Ra N$$ as a Lie subalgebroid of 
$A\times TN$.  Projection to the second factor is a Lie algebroid morphism $f^!A\to A$. 
Under composition, $(f\circ g)^!A=g^! (f^!A)$,  provided that the transversality conditions are satisfied. 
	\begin{enumerate}
		\item For tangent bundles, $f^!TM=TN.$
		 \item  For  Atiyah algebroids,  
		  \[ f^!\on{At}(P)=\on{At}(f^*P)\]
		  with the usual pullback of a principal bundle. 
		\item If $f$ is the inclusion of a submanifold $N\subset M$, then   $f^!A=\a^{-1}(TN)$.
		\item If $f$ is a local diffeomorphism, then $f^!A$ is the pullback as a vector bundle. 
		\item 
		If $\pr_1\colon M\times F\to  M$ is projection to the first factor, then $\pr_1^!A=A\times TF$. This gives the local description of pullbacks under submersions.
	\end{enumerate}

\subsection{Other constructions} 
There are various interesting ways of producing new Lie algebroids out of given ones. The tangent bundle of 
a Lie algebroid $A\Ra M$ has a  Lie algebroid structure $TA\Ra TM$, in such a way that the map on sections $\sigma\mapsto T\sigma$ is a Lie algebra morphism. Similarly,  the $k$-th jet bundle has a Lie algebroid structure $J^k(A)\Ra M$, in such a way that $\sigma\mapsto j^k(\sigma)$ 
is a Lie algebra morphism \cite{kum:lie}.  Given a Lie subalgebroid $B\subset A$  along a submanifold $N\subset M$, the normal bundle is a Lie algebroid 
\[ \nu(A,B)\Ra \nu(M,N).\]
If $N$ is closed, 
there is a \emph{blow-up Lie algebroid} \cite{deb:blo}
\[ \on{Bl}_0(A,B)\Rightarrow \on{Bl}(M,N)\]
replacing the subbundle $A|_N\to N$ with the projectivized normal bundle 
\[ \mathbb{P}(\nu(A,B))\Ra \mathbb{P}(\nu(M,N)).\]
(The subscript indicates that this is different from $\on{Bl}(A,B)$, where one would only remove $B$.) 
 The Lie algebroid structure is such that the blow-down map is a Lie algebroid morphism. If $N=Z$ is a hypersurface, then $\on{Bl}(M,Z)=M$, and so $\on{Bl}_0(A,B)$ is a Lie algebroid over $M$ itself. Special cases include 
$\bb TM$ and $\zz TM$ (see \ref{subsec:defex}\eqref{it:f}) as blow-ups of $TM$ along $TZ$ and $0_Z$, respectively.

\section{Orbits, splitting theorem}
\subsection{Orbits of a Lie algebroid}
Following \cite{and:hol}, we define a \emph{singular foliation} on $M$ to be a subsheaf $\E\subset \mf{X}_M$ of the sheaf of vector fields, 
such that $\E$ is locally finitely generated and involutive: 
$$[\E,\E]\subset \E.$$ 
The first condition means that every point in $M$ has a neighborhood $U$ such that $\E(U)$ is spanned, as a $C^\infty(U)$-module, by finitely many vector fields in $\E(U)$. A \emph{leaf} of a singular foliation is a maximal connected submanifold  $\O\subset M$ whose tangent bundle is 
spanned by vector fields in $\E$; it is a theorem  \cite{and:hol} that every singular foliation defines a decomposition of $M$ into leaves.  
A singular foliation is a \emph{regular} foliation if the leaves have constant dimension. 

Every Lie algebroid $A\Ra M$ determines a singular foliation, with $\E(U)=\a(\Gamma(A|_U))$. Its leaves are  called the \emph{orbits} of the Lie algebroid $A$. The tangent space to the orbit through $x$ is $T_x\O=A|_x/\h_x$ where 
\[ \h_x=\ker(\a_x)\Ra \{x\}\]
is the \emph{isotropy Lie algebra} at $x\in M$. 
It is presently unknown\cite{lau:uni} whether every singular foliation is locally the orbit foliation of a Lie algebroid. 
\begin{enumerate}
	\item The orbits of an action Lie algebroid $\g\times M\Ra M$  are the usual orbits of the $\g$-action.
	\item The orbits of a foliation Lie algebroid  are the leaves of the foliation. 
	\item
	For a Poisson manifold $M$, the orbits of $T^*M\Ra M$ are the symplectic leaves of the Poisson structure. 
	\item The orbits of an Atiyah algebroid $\on{At}(P)\Ra M$ are the components of $M$. 
	\item The orbits of $\bb TM=\on{Bl}_0(TM,TZ)$ are the components of $M-Z,\ Z$. 
\end{enumerate}

A Lie algebroid is called \emph{regular} if it defines a regular  foliation. 
Equivalently, this means that the anchor has constant rank, and so the isotropy Lie algebras 
fit into a family of Lie algebras $\h=\ker(\a)\subset A$.  A Lie algebroid is \emph{transitive} if the anchor is 
surjective, so all its orbits are open. Examples include the action algebroids of transitive $\g$-actions, the Atiyah algebroids of principal bundles,  and the Almeida-Molino algebroid of a closed 2-form.

\subsection{Splitting theorem}
Given a Lie algebroid $A\Ra M$ a submanifold $i\colon N\to M$ is called a \emph{Lie algebroid tranversal}  if it is 
transverse to the anchor (equivalently, transverse to all orbits). In this case, letting $p\colon \nu(M,N)\to N$ be the normal bundle, the Lie algebroid $p^!i^!A\Ra \nu(M,N)$  is the \emph{linear approximation} of $A$ along $N$. 
 
 \begin{theorem}[Linearization along transversals]\cite{bur:spl} 
 There exists an open neighborhood $U\subset \nu(M,N)$ of the zero section and a tubular neighborhood  embedding 
 $\phi\colon U \to M$, lifting to a Lie algebroid isomorphism 
 \[ (p^!i^!A)|_U\to A|_{\phi(U)}.\]	
 \end{theorem}
 The following special case is particularly important:
  
\begin{theorem}[Splitting theorem]\cite{duf:nor,fer:lie, wei:alm} 
Let $\O\subset M$ be an orbit of the Lie algebroid $A\Ra M$, and $N\subset M$ a transversal of complementary dimension, 
with $N\cap \O=\{x\}$. 
Then $A$ is isomorphic, near $x$, to the 
 product Lie algebroid  $$i^!A\times T\O\Ra N\times \O.$$ 
\end{theorem}

If $A\Ra M$ is a regular Lie algebroid, then $i^!A\Ra N$  has zero 
anchor, and so is given by a family of Lie algebras. In particular, if $A\Ra M$ is transitive, then $N=\{x\}$ is just a point, and so we find that $A$ is locally a product 
	$$A|_U\cong \k\times TU$$
with $\k=\h_x$. In particular, 	$\ker(\a)\subset A$ is a Lie algebra bundle, with typical fiber $\k$.

\section{The correspondence between Lie groupoids and Lie algebroids}
\subsection{The Lie functor}
%
Let $G\rra M$ be a Lie groupoid, with source and target maps denoted
$\sz,\tz\colon G\to M$. A vector field $X\in\mf{X}(G)$ is \emph{left-invariant} if it is tangent to the $\tz$-fibers, and 
$$ X_g=(L_g)_* X_{\sz(g)}$$
for all $g\in G$, where $L_g\colon  \tz^{-1}(\sz(g))\mapsto \tz^{-1}(\tz(g))$ denotes left translation. The Lie algebroid 
$$A=\on{Lie}(G)\Ra M$$ 
of a Lie groupoid is defined as follows:

\begin{itemize}
	\item As a vector bundle, $A=\ker(T\tz)|_M$.
	\item The bracket on $\Gamma(A)$ is defined by its identification  with the left-invariant vector fields. 
	\item The anchor is the restriction of $-T\sz\colon TG\to TM$ to $A\subset TG|_M$.  
\end{itemize} 
Some examples:
\begin{enumerate}
	\item For a Lie group one recovers the usual construction of the Lie algebra. 
	\item The pair groupoid  and the homotopy groupoid of $M$ both have Lie algebroid $TM\Ra M$.
	\item The Lie algebroid of the holonomy groupoid or monodromy groupoid of a foliation $\ca{F}$ is $T_\F M\Ra M$.  
	\item The Lie algebroid of the gauge groupoid  of a principal bundle $P\to M$ is the Atiyah algebroid $\on{At}(P)\Ra M$.  
\end{enumerate}

\subsection{Integrability}
A Lie algebroid $A\Ra M$ is called \emph{integrable} if it is of the form $A=\on{Lie}(G)$ for a Lie groupoid $G\rra M$. 
Examples from foliation theory dictate that one should  allow for Lie groupoids $G\rra M$  whose total space is \emph{non-Hausdorff} (but still insisting that the manifold of units and all $\sz,\tz$-fibers are Hausdorff. If $A$ is integrable, then there is a unique integration whose source fibers are simply connected. 
\begin{enumerate}
	\item The source-simply connected integration of $TM\Ra M$  is the homotopy groupoid $\pi(M)\rra M$.  
	\item Given a foliation $\F$, the source-simply connected integration of $T_\F M\Ra M$  is the \emph{holonomy groupoid}. \cite{pra:feu}, \cite{win:gra}. 
	\item Every family of Lie algebras (as a Lie algebroid with zero anchor) is integrable to a family of Lie groups \cite{dou:esp}.
	\item Every action Lie algebroid $\k\times M\Ra M$ is integrable (but the integration is not an action groupoid for a global group action, in general) \cite{daz:int}.  
\end{enumerate}

A counter-example to integrability was observed by Almeida-Molino \cite{alm:sui1}, who found that the Lie algebroid 
$A_\omega$ defined by a closed 2-form $\omega$ (see \ref{subsec:defex}\eqref{it:g}) 
cannot be integrable unless the group of spherical periods 
\[ \Gamma=\{\int_{S^2}j^*\omega|\ j\colon S^2\to M\}\subset \R\]
is discrete. 
The precise criterion for integrability of Lie algebroids $A\Ra M$ will be given below. 
It involves the \emph{monodromy invariants}. 

\subsection{Monodromy groups}
The original definition \cite{cra:intlie} of the monodromy groups uses Lie algebroid homotopy. Here we present an alternative approach
\cite{me:int}. 
Let $K$ be a connected Lie group. Recall that isomorphism classes of principal bundles over 2-spheres $S^2$, with trivialization $P_{x_0}\cong K$ at some fixed base point $x_0$, are labeled by the fundamental group $\pi_1(K)$. This is a group isomorphism, for the group structure on classes of principal bundles given by connected sum. The proof is by a clutching construction, and may be used 
to obtain  a similar result for Lie algebroids. Let $\k$ be a Lie algebra,  denote by $\wt{K}$ the connected and simply connected 
Lie group  integrating $\k$, and let $\on{Cent}(\wt{K})$ be its center. 
\begin{theorem}\cite{me:int}\label{th:me}
	There is a canonical isomorphism between  the group of equivalence classes of transitive Lie algebroids $A\Ra S^2$, together 
	with an isomorphism $\ker(\a_{x_0})\cong \k$, and the group $\on{Cent}(\wt{K})$. 
	Given a connected Lie group $K$ with Lie algebra $\k$, the natural map 
	$$ \pi_1(K)\to \on{Cent}(\wt{K})$$ 
	takes the class of a principal $K$-bundle to the class of its Atiyah algebroid. 
\end{theorem}
For instance, since $\on{Cent}(\SU(2))=\Z_2$, there are two classes of  Lie algebroids with 
$\k=\mf{su}(2)$ .  They are represented by the Atiyah algebroids of the trivial and non-trivial principal $\SO(3)$-bundles.
On the other hand, for $\k=\R$ we obtain a 1-parameter family of classes of Lie algebroid, represented by the Almeida-Molino Lie algebroids for multiples of the volume form.

Suppose $A\Ra M$ is any transitive Lie algebroid, and $x\in M$, with isotropy Lie algebra  $\h_x=\ker(\a_x)$. 
Given a smooth map $j\colon (S^2,x_0)\to (M,x)$, the pullback 
$j^!A\Ra S^2$ is transitive, and so defines an element of $\on{Cent}(\wt{H}_x)$. 
The set of elements obtained in this way define the  \emph{monodromy group}  
\[ \Gamma_x\subset \on{Cent}(\wt{H}_x).\]
By transitivity, the monodromy  groups for different points $x\in M$ are isomorphic. For general, not necessarily transitive Lie algebroids $A\Ra M$, we may define the monodromy group at $x\in M$ as the 
monodromy group $\Gamma_x$ of the restriction $A_\O$ to the orbit $\O$ through $x$.

\subsection{Crainic-Fernandes integrability theorem}
Consider first the case of a transitive Lie algebroid $A\Ra M$.  Fix $x\in M$, with isotropy algebra $\h_x$. If $A$ is integrable, then its integration is a transitive Lie groupoid $G\rra M$, and as such is the gauge groupoid 
of a principal $K$-bundle $P\to M$, where $K$ is a Lie group (possibly disconnected) with Lie algebra $\k=\h_x$. Accordingly, 
$A$ is the Atiyah algebroid $\on{At}(P)$.   By the 
second part of \ref{th:me}, and since $j^!\on{At}(P)=\on{At}(j^*P)$, the elements of the monodromy groups are all contained in 
$\pi_1(K^0)\subset \on{Cent}(\wt{K})$. In particular, the monodromy group $\Gamma_x$ must be discrete. 

It turns out that this condition is also sufficient. In fact, there is an explicit construction \cite{me:int} of a source-simply connected Lie groupoid integrating $A$, analogous to the description of a principal bundle in terms of the parallel transport with respect to a connection. 

For non-transitive Lie algebroids $A\Ra M$ discreteness of the monodromy groups is a necessary condition for integrability, since an integration of $A$ also gives an integration of its restriction to orbits. However, it is not sufficient, in general.  
%
%
\begin{theorem}[Crainic-Fernandes] \cite{cra:intlie}\label{th:crafe}
A Lie algebroid $A\Ra M$ is integrable if and only if the monodromy groups are \emph{uniformly discrete}, in the sense that 
there is an open neighborhood $U$ of the zero section of $A$ with the property that if $\xi\in \h_x\subset A_x$ is contained 
in $U$, 
and $\exp\xi\in \Gamma_x$, then $\xi=0$. 
\end{theorem}

Again, when these conditions are satisfied, there is an explicit construction of the source-simply connected 
groupoid integrating $A$, as a quotient 
\[ G=\frac{\{\mbox{Lie algebroid paths}\ \ T[0,1]\to A\}}{\{\mbox{Lie algebroid homotopy}\}}.\]
A groupoid structure is defined by concatenation. The construction was suggested by Catteneo-Felder \cite{cat:poi} in the case of cotangent Lie algebroids
and by \v{S}evera \cite{sev:so} and Weinstein for general Lie algebroids. It is motivated by 
the proof of Lie's third theorem for Lie groups by Duistermaat-Kolk \cite{du:li}. As shown in   \cite{cra:intlie}, 
the  groupoid $G$ is always defined as a topological groupoid; Theorem \ref{th:crafe} gives a criterion for its smoothness.

\section{The Lie algebroid complex}
\subsection{Contravariant descriptions of Lie algebroids}
Given a vector bundle $A\to M$, consider the graded algebra 
$\mathsf{C}^\bullet(A)=\bigoplus_k \mathsf{C}^k(A)$ with 
\[ \mathsf{C}^k(A)=\Gamma(\wedge^k A^*).\]

%
\begin{theorem}\label{th:covariant} \cite{vai:lie}
 	There is a bijective correspondence between Lie algebroid structures on  vector bundles $A$ and degree $1$ 
		differentials on the graded algebra $C^\bullet(A)$. 
		A vector bundle morphism $B\to A$ between Lie algebroids is a Lie algebroid morphism if and only if 
		the induced map $\mathsf{C}^\bullet(A)\to \mathsf{C}^\bullet(B)$ is a cochain map.  
\end{theorem}
In the language of super-geometry,  $\mathsf{C}^\bullet(A)$ is the algebra of smooth functions on the super-manifold $A[1]$, and the differential is a homological vector field $Q$ on $A[1]$.  Given a Lie algebroid structure on $A$, the differential is given by 
%
\begin{align*} (\d\phi)(\sigma_1,\ldots,\sigma_{k+1})
&=
\sum_i (-1)^{i+1}\L_{\a(\sigma_i)}\phi(\sigma_1,\ldots,\wh{\sigma_i},\ldots)\\
&+\sum_{i<j} (-1)^{i+j}\phi([\sigma_i,\sigma_j],\sigma_1,\ldots,\wh{\sigma_i},\ldots,\wh{\sigma_j},\ldots).
\end{align*}
The complex $(\mathsf{C}^\bullet(A),\d)$
generalizes both the de Rham complex $\Omega(M)$ of differential forms (for $A=TM$) and the Chevalley-Eilenberg complex $\wedge\g^*$ for $A=\g$. If $A=T_\F M$ for a foliation $\F$, one obtains the leaf-wise de Rham complex. 
If $M$ is a Poisson manifold, the Lie algebroid complex for $T^*M\Ra M$
is given by the Lichnerowicz complex $\mathsf{C}^\bullet(T^*M)=\mf{X}^\bullet(M)$, defining \emph{Poisson cohomology}. 
As another formulation of Theorem \ref{th:covariant}, we have:
\begin{theorem}
	Lie algebroid structures on vector bundles $A\to M$ are equivalent to Poisson 
	structures on the dual bundle $A^*$ which are \emph{linear}, in the sense that  the Poisson bracket between linear functions on $A$ is again linear.
\end{theorem} 
Identifying on $M$ with their pullback to $A^*$, and letting $\phi_\sigma\in C^\infty(A^*)$ be the linear function defined 
by the section $\sigma$, the Poisson bracket is determined by 
\[ \{f,g\}=0,\ \ \{\phi_\sigma,f\}=\L_{\a(\sigma)}f,\ \ \{\phi_\sigma,\phi_\tau\}=\phi_{[\sigma,\tau]}.\]
In particular, for the cotangent Lie algebroid $A=T^*M$ of a Poisson manifold $M$ one obtains a linear Poisson structure on $A^*=TM$ called the tangent lift of the Poisson structure. 	

\subsection{Van Est map}
For Lie groupoids $G\rra M$, there is a complex $\mathsf{C}^\bullet(G)$ of normalized groupoid cochains, 
generalizing the complex of smooth Lie group cochains. If $A=\on{Lie}(G)$, there is a \emph{Van Est map},
\[ \on{VE}^\bullet \colon \mathsf{C}^\bullet(G)\to \mathsf{C}^\bullet(A),\]
due to Weinstein-Xu \cite{wei:ext}, which is a morphism of differential graded algebras. 
As shown by Crainic \cite{cra:dif}, in generalization of the van Est theorem for Lie groups, the induced map in cohomology is an isomorphism up to degree $k$ and injective in degree $k+1$, provided that the  source fibers of $G$ are $k$-connected. Van Est  integration formulas were obtained in \cite{cab:loc,mei:ve2}.

\section{Outlook}
Due to limitations of space, we had to leave out many relevant aspects of Lie algebroids. Here are several such topics: Universal enveloping algebras of Lie algebroids \cite{nis:pse}, Lie algebroid connections and characteristic classes \cite{fer:lie}, modular classes \cite{ev:tra}, geometric structures on Lie algebroids and their integration  \cite{bur:mul}, \cite{mac:lif}, Lie bialgebroids and their integration \cite{kos:exa,mac:int}, 
$\VB$-algebroids \cite{bur:vec,gra:vba} and double Lie algebroids \cite{mac:dou,vor:q}, Dirac structures and Courant algebroids \cite{couwein:beyond,lib:cou},  applications in generalized complex geometry \cite{gua:ge1}, local integration of Lie algebroids \cite{cab:loc}, 
Lie-$n$ algebroids and shifted Lie algebroids \cite{lau:uni,pym:sh,she:hig}. As a general reference on Lie algebroids we suggest the book of Mackenzie \cite{mac:gen}.
\bibliographystyle{amsplain} 

\begin{thebibliography}{10}
 	
 	\bibitem{alm:sui1}
 	R.~Almeida and P.~Molino, \emph{{Suites d'Atiyah, feuilletages et
 			quantification g\'eom\'etrique. (Atiyah sequences, foliations and geometric
 			quantification)}}, {S\'emin. G\'eom. Diff\'er., Univ. Sci. Tech. Languedoc
 		1984/1985, 39-59 (1985).}, 1985.
 	
 	\bibitem{and:hol}
 	I.~Androulidakis and G.~Skandalis, \emph{The holonomy groupoid of a singular
 		foliation}, J. Reine Angew. Math. \textbf{626} (2009), 1--37.
 	
 	\bibitem{bur:mul}
 	H.~Bursztyn and A.~Cabrera, \emph{Multiplicative forms at the infinitesimal
 		level}, Math. Ann. \textbf{353} (2012), no.~3, 663--705.
 	
 	\bibitem{bur:vec}
 	H.~Bursztyn, A.~Cabrera, and M.~del Hoyo, \emph{{Vector bundles over Lie
 			groupoids and algebroids}}, Adv.~Math. \textbf{290} (2016), no.~2, 163--207.
 	
 	\bibitem{bur:spl}
 	H.~Bursztyn, H.~Lima, and E.~Meinrenken, \emph{{Splitting theorems for Poisson
 			and related structures}}, J. Reine Angew.~Math.~ \textbf{754} (2019),
 	281--312.
 	
 	\bibitem{cab:loc}
 	A.~Cabrera, I.~Marcut, and A.~Salazar, \emph{{On local integration of Lie
 			brackets}}, Journal f\"ur die reine und angewandte Mathematik \textbf{760}
 	(2020), 267--293.
 	
 	\bibitem{cat:poi}
 	A.~Cattaneo and G.~Felder, \emph{Poisson sigma models and symplectic
 		groupoids}, Quantization of singular symplectic quotients, Progr. Math., vol.
 	198, Birkh\"auser, Basel, 2001, pp.~61--93.
 	
 	\bibitem{couwein:beyond}
 	T.~Courant and A.~Weinstein, \emph{Beyond {P}oisson structures}, Action
 	hamiltoniennes de groupes.~Troisi\`eme th\'eor\`eme de {L}ie (Lyon, 1986),
 	Travaux en Cours, vol.~27, Hermann, Paris, 1988, pp.~39--49.
 	
 	\bibitem{cra:dif}
 	M.~Crainic, \emph{Differentiable and algebroid cohomology, van {E}st
 		isomorphisms, and characteristic classes}, Comment. Math. Helv. \textbf{78}
 	(2003), no.~4, 681--721.
 	
 	\bibitem{cra:intlie}
 	M.~Crainic and R.L. Fernandes, \emph{Integrability of {L}ie brackets}, Ann. of
 	Math. (2) \textbf{157} (2003), no.~2, 575--620.
 	
 	\bibitem{daz:int}
 	P.~Dazord, \emph{Int\'egration d'alg\`ebres de {L}ie locales et groupo\"\i des
 		de contact}, C. R. Acad. Sci. Paris S\'er. I Math. \textbf{320} (1995),
 	no.~8, 959--964.
 	
 	\bibitem{deb:blo}
 	C.~Debord and G.~Skandalis, \emph{{Blowup constructions for Lie groupoids and a
 			Boutet de Monvel type calculus}}, M\"unster Journal of Mathematics
 	\textbf{14} (2021), 1--40.
 	
 	\bibitem{dou:esp}
 	A.~Douady and M.~Lazard, \emph{Espaces fibr\'es en alg\`ebres de {L}ie et en
 		groupes}, Invent. Math. \textbf{1} (1966), 133--151.
 	
 	\bibitem{duf:nor}
 	J.-P. Dufour, \emph{Normal forms for {L}ie algebroids}, Lie algebroids and
 	related topics in differential geometry ({W}arsaw, 2000), Banach Center
 	Publ., vol.~54, Polish Acad. Sci. Inst. Math., Warsaw, 2001, pp.~35--41.
 	
 	\bibitem{du:li}
 	J.~J. Duistermaat and J.~A.~C. Kolk, \emph{{Lie Groups}}, Springer-Verlag,
 	Berlin, 2000.
 	
 	\bibitem{ev:tra}
 	S.~Evens, J.-H. Lu, and A.~Weinstein, \emph{Transverse measures, the modular
 		class and a cohomology pairing for {L}ie algebroids}, Quart.~J.~Math.~Oxford
 	Ser.~(2) \textbf{50} (1999), no.~200, 417--436.
 	
 	\bibitem{fer:lie}
 	R.~Fernandes, \emph{Lie algebroids, holonomy and characteristic classes}, Adv.
 	Math. \textbf{170} (2002), no.~1, 119--179.
 	
 	\bibitem{gra:vba}
 	A.~Gracia-Saz and R.~Mehta, \emph{{Lie algebroid structures on double vector
 			bundles and representation theory of Lie algebroids}}, Advances in
 	Mathematics \textbf{223} (2010), 1236--1275.
 	
 	\bibitem{gua:ge1}
 	M.~Gualtieri, \emph{Generalized complex geometry}, Ann. of Math. (2)
 	\textbf{174} (2011), no.~1, 75--123.
 	
 	\bibitem{kos:exa}
 	Y.~Kosmann-Schwarzbach, \emph{Exact {G}erstenhaber algebras and {L}ie
 		bialgebroids}, Acta Appl. Math. \textbf{41} (1995), no.~1-3, 153--165,
 	Geometric and algebraic structures in differential equations.
 	
 	\bibitem{kum:lie}
 	A.~Kumpera and D.~Spencer, \emph{Lie equations. {V}ol. {I}: {G}eneral theory},
 	vol. No. 73., Princeton University Press, Princeton, N.J.; University of
 	Tokyo Press, Tokyo, 1972.
 	
 	\bibitem{lau:uni}
 	C.~Laurent-Gengoux, S.~Lavau, and T.~Strobl, \emph{The universal {L}ie
 		{$\infty$}-algebroid of a singular foliation}, Doc. Math. \textbf{25} (2020),
 	1571--1652.
 	
 	\bibitem{lib:cou}
 	D.~Li-Bland and E.~Meinrenken, \emph{{C}ourant algebroids and {P}oisson
 		geometry}, International Mathematics Research Notices \textbf{11} (2009),
 	2106--2145.
 	
 	\bibitem{mac:dou}
 	K.~Mackenzie, \emph{Double {L}ie algebroids and second-order geometry.~{I}},
 	Adv.~Math. \textbf{94} (1992), no.~2, 180--239.
 	
 	\bibitem{mac:gen}
 	\bysame, \emph{General theory of {L}ie groupoids and {L}ie algebroids}, London
 	Mathematical Society Lecture Note Series, vol. 213, Cambridge University
 	Press, Cambridge, 2005.
 	
 	\bibitem{mac:lif}
 	K.~Mackenzie and P.~Xu, \emph{Classical lifting processes and multiplicative
 		vector fields}, Quart. J. Math. Oxford Ser. (2) \textbf{49} (1998), no.~193,
 	59--85.
 	
 	\bibitem{mac:int}
 	\bysame, \emph{Integration of {L}ie bialgebroids}, Topology \textbf{39} (2000),
 	no.~3, 445--467.
 	
 	\bibitem{me:int}
 	E.~Meinrenken, \emph{On the integration of transitive {L}ie algebroids},
 	Enseign. Math. \textbf{67} (2021), no.~3-4, 423--454.
 	
 	\bibitem{mei:ve2}
 	E.~Meinrenken and M.~A. Salazar, \emph{Van {E}st differentiation and
 		integration}, Math. Ann. \textbf{376} (2020), no.~3-4, 1395--1428.
 	\MR{4081118}
 	
 	\bibitem{mel:ati}
 	R.~B. Melrose, \emph{{The Atiyah-Patodi-Singer index theorem}}, Research Notes
 	in Mathematics, vol.~4, A K Peters, Ltd., Wellesley, 1993.
 	
 	\bibitem{nis:pse}
 	V.~Nistor, A.~Weinstein, and P.~Xu, \emph{Pseudodifferential operators on
 		differential groupoids}, Pacific J. Math. \textbf{189} (1999), no.~1,
 	117--152.
 	
 	\bibitem{pra:th}
 	J.~Pradines, \emph{Th\'eorie de {L}ie pour les groupo\"\i des
 		diff\'erentiables. {C}alcul diff\'erenetiel dans la cat\'egorie des
 		groupo\"\i des infinit\'esimaux}, C. R. Acad. Sci. Paris S\'er. A-B
 	\textbf{264} (1967), A245--A248.
 	
 	\bibitem{pra:feu}
 	\bysame, \emph{Feuilletages: holonomie et graphes locaux}, C. R. Acad. Sci.
 	Paris S\'{e}r. I Math. \textbf{298} (1984), no.~13, 297--300. \MR{765427}
 	
 	\bibitem{pym:sh}
 	B.~Pym and P.~Safronov, \emph{{Shifted symplectic Lie algebroids}}, Int. Math.
 	Res. Not. IMRN (2020), no.~21, 7489--7557.
 	
 	\bibitem{sev:so}
 	P.~{\v{S}}evera, \emph{{Some title containing the words `homotopy' and
 			`symplectic', e.g. this one}}, Travaux Math\'ematiques, Univ. Luxemb.
 	\textbf{XVI} (2005), 121--137, based on talk at Poisson 2000, Luminy.
 	
 	\bibitem{she:hig}
 	Y.~Sheng and C.~Zhu, \emph{Higher extensions of {L}ie algebroids}, Commun.
 	Contemp. Math. \textbf{19} (2017), no.~3, 1650034, 41.
 	
 	\bibitem{vai:lie}
 	A.~Vaintrob, \emph{Lie algebroids and homological vector fields}, Uspekhi Mat.
 	Nauk \textbf{52} (1997), no.~2(314), 161--162.
 	
 	\bibitem{vor:q}
 	T.~Voronov, \emph{{$Q$}-manifolds and {M}ackenzie theory}, Comm. Math. Phys.
 	\textbf{315} (2012), no.~2, 279--310.
 	
 	\bibitem{wei:alm}
 	A.~Weinstein, \emph{Almost invariant submanifolds for compact group actions},
 	J. Eur. Math. Soc. (JEMS) \textbf{2} (2000), no.~1, 53--86.
 	
 	\bibitem{wei:ext}
 	A.~Weinstein and P.~Xu, \emph{Extensions of symplectic groupoids and
 		quantization}, J. Reine Angew. Math. \textbf{417} (1991), 159--189.
 	
 	\bibitem{win:gra}
 	H.~E. Winkelnkemper, \emph{The graph of a foliation}, Ann. Global Anal. Geom.
 	\textbf{1} (1983), no.~3, 51--75. \MR{739904}
 	
 \end{thebibliography}
 
 \def\cprime{$'$} \def\polhk#1{\setbox0=\hbox{#1}{\ooalign{\hidewidth
 			\lower1.5ex\hbox{`}\hidewidth\crcr\unhbox0}}} \def\cprime{$'$}
 \def\cprime{$'$} \def\cprime{$'$} \def\cprime{$'$} \def\cprime{$'$}
 \def\polhk#1{\setbox0=\hbox{#1}{\ooalign{\hidewidth
 			\lower1.5ex\hbox{`}\hidewidth\crcr\unhbox0}}} \def\cprime{$'$}
 \def\cprime{$'$} \def\cprime{$'$} \def\cprime{$'$} \def\cprime{$'$}
 \providecommand{\bysame}{\leavevmode\hbox to3em{\hrulefill}\thinspace}
 \providecommand{\MR}{\relax\ifhmode\unskip\space\fi MR }
 \providecommand{\MRhref}[2]{%
 	\href{http://www.ams.org/mathscinet-getitem?mr=#1}{#2}
 }
 \providecommand{\href}[2]{#2}

\end{document}